\documentstyle[12pt]{article}
\font\teneufm=eufm10
\font\seveneufm=eufm7
\font\fiveeufm=eufm5
\newfam\eufmfam
\textfont\eufmfam=\teneufm
\scriptfont\eufmfam=\seveneufm
\scriptscriptfont\eufmfam=\fiveeufm
\def\frak#1{{\fam\eufmfam\relax#1}}

\newfam\msbfam
\font\tenmsb=msbm10 scaled \magstep1  \textfont\msbfam=\tenmsb
\font\sevenmsb=msbm7 scaled \magstep1 \scriptfont\msbfam=\sevenmsb
\font\fivemsb=msbm5 scaled \magstep1  \scriptscriptfont\msbfam=\fivemsb
\def\Bbb{\fam\msbfam \tenmsb}

\def\RR{{\Bbb R}}
\def\CC{{\Bbb C}}

\def\Im{\hbox{Im}\,}
\def\ra{\rightarrow}

\def\HollowBoxx #1#2#3{{\dimen0=#1 \advance\dimen0 by -#2       
       \dimen1=#1 \advance\dimen1 by #3                       
        \vrule height 0pt depth #3 width #2                   
       \hskip -#3
       \vrule height #1 depth #3 width #3}}                   
 \def\LeftContraction{\mathord{\kern1.45pt \HollowBoxx{6pt}{3.5pt}{.4pt}}\,}

 \def\HollowBox #1#2#3{{\dimen0=#1 \advance\dimen0 by -#3       
       \dimen1=#1 \advance\dimen1 by #3                       
        \vrule height #1 depth #3 width #3                    
        \vrule height 0pt depth #3 width #2                   
        \hskip -#3}}                                             
 \def\RightContraction{\mathord{\, \HollowBox{6pt}{3.1pt}{.4pt}}\kern1.6pt}

\newtheorem{theorem}{THEOREM}[section]

\newtheorem{lemma}[theorem]{Lemma}

\begin{document}
\begin{center}
{\Large \bf Canonical Isomorphism of Two Lie Algebras} 
\medskip \\
{\Large \bf Arising in $CR$-geometry}\footnote{{\bf Mathematics 
    Subject Classification:} 32C16, 32F40, 17B66}\footnote{{\bf Keywords 
   and Phrases:} $CR$-manifolds, equivalence problem, isomorphism of
    Lie algebras}  
\medskip \\
{\normalsize V. V. Ezhov and A. V. Isaev}
\end{center} 

\begin{quotation} 
\small \sl We show that the maximal prolongation of a certain
algebra associated with a non-degenerate Hermitian form on $\CC^n\times\CC^n$
with values in $\RR^k$ is canonically isomorphic to the Lie algebra of
infinitesimal holomorphic automorphisms of the corresponding quadric in $\CC^{n+k}$.
This fact creates a link between different approaches to the
equivalence problem for Levi-nondegenerate strongly uniform $CR$-manifolds. 
\end{quotation}

\markboth{V. V. Ezhov and A. V. Isaev}{Canonical Isomorphism of Two Lie Algebras}

\setcounter{section}{-1}

\section{Introduction and Formulation of Result}

A $CR$-{\it structure} on a smooth real manifold $M$ of dimension $m$ is a smooth distribution of
subspaces in the tangent spaces $T_p^c(M)\subset T_p(M)$, $p\in M$,
with operators of complex structure $J_p:T_p^c(M)\ra T_p^c(M)$,
$J_p^2\equiv -\hbox{id}$, that depend smoothly on $p$. A manifold $M$
equipped with a $CR$-structure is called a $CR$-{\it manifold}. It
follows that the number
$CR\hbox{dim}M:=\hbox{dim}_{\CC}T_p^c(M)$ does not depend on
$p$; it is called the $CR$-{\it dimension} of  $M$. The number
$CR\hbox{codim}M:=m-2CR\hbox{dim}M$ is called the $CR$-{\it
codimension} of $M$. $CR$-structures naturally arise on real submanifolds in complex manifolds. Indeed, if, for example,
$M$ is a real submanifold of $\CC^K$, then one can define the distribution
$T_p^c(M)$ as follows:
$$
T_p^c(M):=T_p(M)\cap iT_p(M).
$$
On each $T_p^c(M)$ the operator $J_p$ is then defined as the
operator of multiplication by $i$. Then $\{T_p^c(M), J_p\}_{p\in M}$
form a $CR$-structure on $M$, if $\dim_{\CC}T_p^c(M)$ is constant. This is always the case, for example, if $M$ is a real
hypersurface in $\CC^K$ (in which case $CR\hbox{codim}M=1$). We say that such a $CR$-structure is {\it
induced by $\CC^K$}.

A mapping between two $CR$-manifolds $f:M_1\ra M_2$ is called a
$CR$-{\it mapping}, if for every $p\in M_1$: {\bf (i)} $df(p)$ maps
$T_p^c(M_1)$ to $T_{f(p)}^c(M_2)$, and {\bf (ii)} $df(p)$ is complex linear
on $T_p^c(M_1)$. Two $CR$-manifolds $M_1$, $M_2$ are called $CR$-{\it
equivalent}, if there is a $CR$-diffeomorphism from $M_1$ onto
$M_2$. Such a $CR$-diffeomorphism $f$ is called a $CR$-{\it
isomorphism}.

Let $M$ be a $CR$-manifold. For
every $p\in M$ consider the complexification
$T_p^c(M)\otimes_{\RR}\CC$. Clearly, this complexification can be
represented as the direct sum
$$
T_p^c(M)\otimes_{\RR}\CC=T_p^{(1,0)}(M)\oplus T_p^{(0,1)}(M),
$$
where
\begin{eqnarray*}
T_p^{(1,0)}(M)&:=&\{X-iJ_pX:X\in T_p^c(M)\},\\
T_p^{(0,1)}(M)&:=&\{X+iJ_pX:X\in T_p^c(M)\}.
\end{eqnarray*}
The $CR$-structure on $M$ is called {\it integrable} if for any local
sections $Z,Z'$ of the bundle $T^{(1,0)}(M)$, the vector field $[Z,Z']$
is also a section of $T^{(1,0)}(M)$. It is not difficult to see that
if $M\subset\CC^K$ and the $CR$-structure on $M$ is induced by
$\CC^K$, then it is integrable.

An important characteristic of a $CR$-structure called the {\it Levi
form} comes from taking commutators of local sections of
$T^{(1,0)}(M)$ and $T^{(0,1)}(M)$. Let $p\in M$, $z, z'\in
T_p^{(1,0)}(M)$, and $Z$, $Z'$ be local sections of $T^{(1,0)}(M)$ near
$p$ such that $Z(p)=z$, $Z'(p)=z'$. The Levi form of $M$ at
$p$ is the Hermitian form on $T_p^{(1,0)}(M)\times T_p^{(1,0)}$ with values in $(T_p(M)/T_p^c(M))\otimes_{\RR}\CC$ given by
$$
{\cal L}_M(p)(z,z'):=i[Z,\overline{Z'}](p)(\hbox{mod}\,T_p^c(M)\otimes_{\RR}\CC).
$$
The Levi form is defined uniquely up to the choice of coordinates in\linebreak
$(T_p(M)/T_p^c(M))\otimes_{\RR}\CC$, and, for fixed $z$ and $z'$,
its value  does not depend on the choice of $Z$ and $Z'$.

Let $H=(H^1,\dots,H^k)$ be a Hermitian form on $\CC^n\times\CC^n$ with values in
$\RR^k$. For any such $H$ there is a corresponding standard $CR$-manifold
$Q_H\subset\CC^{n+k}$ of $CR$-dimension $n$ and $CR$-codimension
$k$ defined as follows:
$$
Q_H:=\{(z,w):\Im w=H(z,z)\},
$$
where $z:=(z_1,\dots,z_n)$, $w:=(w_1,\dots,w_k)$ are coordinates in
$\CC^{n+k}$.
The manifold $Q_H$ is often called the {\it quadric associated with
the form $H$}. The Levi form of $Q_H$ at any point is given by
$H$.

A Hermitian form $H$ is called {\it non-degenerate} if:
\smallskip\\

\noindent{\bf (i)} The scalar Hermitian forms $H^1,\dots,H^k$ are
linearly independent over $\RR$;
\smallskip\\

\noindent{\bf (ii)} $H(z,z')=0$ for all $z'\in\CC^n$ implies $z=0$.
\smallskip\\

A $CR$-structure on $M$ is
called {\it Levi non-degenerate}, if its
Levi form at any $p\in M$ is non-degenerate.
An important tool in the geometry of Levi non-degenerate integrable $CR$-manifolds is
the {\it automorphism group of $Q_H$}. Let $\hbox{Aut}(Q_H)$ denote
the collection of all local $CR$-isomorphisms of $Q_H$ to itself that
we call {\it local $CR$-automorphisms}. It
turns out that, if $H$ is non-degenerate, then any
local $CR$-automorphism extends to a rational (more precisely, a matrix
fractional quadratic) map of $\CC^{n+k}$ \cite{KT},
\cite{F}, \cite{Tu}, \cite{ES1}. Thus, for a non-degenerate $H$, $\hbox{Aut}(Q_H)$
is a finite-dimensional Lie group. Let ${\frak g}_H$ denote the Lie
algebra of $\hbox{Aut}(Q_H)$. As shown in \cite{B},
\cite{S}, (see also \cite{ES2} for a simple proof), the algebra ${\frak g}_H$ consists of
polynomial vector fields on $\CC^{n+k}$ of the form
\begin{eqnarray*}
{\frak
g}_H&=&\Biggl\{\Biggl(p+Cz+aw+A(z,z)+B(z,w)\Biggr)\frac{\partial}{\partial
z}+\\
&{}&\Biggl(q+2iH(z,p)+sw+2iH(z,a\overline{w})+r(w,w)\Biggr)
\frac{\partial}{\partial w}\Biggr\},\qquad (0.1)
\end{eqnarray*}
where $p\in\CC^n$, $q\in\RR^k$, $C$ is an $n\times n$-matrix, $s$ is a
$k\times k$-matrix, $A(z,z)$ is a quadratic form on $\CC^n\times\CC^n$ with
values in $\CC^n$, $a$ is an $n\times k$-matrix, $B(z,w)$ is a
bilinear form on $\CC^n\times\CC^k$ with values in $\CC^n$, $r(w_1,w_2)$
is a symmetric bilinear form on $\CC^k\times\CC^k$ with values in $\CC^k$, and the
following holds
\begin{eqnarray*}
2\hbox{Re}\,H(Cz,z)&=&sH(z,z),\qquad\qquad\qquad\qquad(0.2.a)\\
H(A(z,z),z)&=&2iH(z,aH(z,z)),\,\,\qquad\qquad(0.2.b)\\
\hbox{Re}\,H(B(z,u),z)&=&r(H(z,z),u),\qquad\qquad\qquad(0.2.c)\\
\hbox{Im}\,H(B(z,H(z,z)),z)&=&0,\,\,\,\,\,\qquad\qquad\qquad\qquad\qquad(0.2.d)
\end{eqnarray*}
for all $z\in\CC^n$, $u\in\RR^k$.

We can now make ${\frak g}_H$ into a graded Lie algebra by introducing
weights as follows: $z$ has weight 1, $w$ has weight 2, $\frac{\partial}{\partial
z}$ has weight -1, $\frac{\partial}{\partial w}$ has weight -2. Then we get ${\frak
g}_H=\oplus_{l=-2}^2{\frak g}^l_H$, where
\begin{eqnarray*}
{\frak g}^{-2}_H&:=&\left\{q\frac{\partial}{\partial
w}\right\},\,\,\,\qquad\qquad\qquad\qquad\qquad\qquad\qquad (0.3.a)\\
{\frak g}^{-1}_H&:=&\left\{p\frac{\partial}{\partial
z}+2iH(z,p)\frac{\partial}{\partial
w}\right\},\qquad\qquad\qquad\qquad (0.3.b)\\
{\frak g}^0_H&:=&\left\{Cz\frac{\partial}{\partial
z}+sw\frac{\partial}{\partial
w}\right\},\qquad\qquad\qquad\qquad\qquad (0.3.c)\\
{\frak g}^1_H&:=&\left\{\left(aw+A(z,z)\right)\frac{\partial}{\partial
z}+2iH(z,a\overline{w})\frac{\partial}{\partial
w}\right\},\,\,\,\,\,\,\,\,(0.3.d)\\
{\frak g}^2_H&:=&\left\{B(z,w)\frac{\partial}{\partial
z}+r(w,w)\frac{\partial}{\partial w}\right\}.\,\,\qquad\qquad\qquad{(0.3.e)}
\end{eqnarray*}
Note that $Q_H$ is a
homogeneous manifold since the global $CR$-automorphisms
\begin{eqnarray*}
z&\mapsto& z+z^0,\\
w&\mapsto&w+w^0+2iH(z,z^0),\qquad\qquad\qquad\qquad {(0.4)}
\end{eqnarray*}
for $(z^0,w^0)\in Q_H$, act transitively on $Q_H$. The subalgebra ${\frak
g}^{-1}_H\oplus{\frak g}^{-2}_H$ is the Lie algebra of the subgroup of
$\hbox{Aut}(Q_H)$ consisting of
automorphisms of the form (0.4). The subalgebra ${\frak g}^0_H$ is the Lie
algebra of the subgroup of $\hbox{Aut}(Q_H)$ consisting of linear automorphisms,
i.e. automorphisms of the form
$$
z\mapsto Pz,\qquad w\mapsto Rw,
$$
where $P$ is a complex $n\times n$-matrix, $R$ is a real $k\times
k$-matrix such that
$$
R^{-1}H(Pz,Pz)=H(z,z).
$$
The components ${\frak g}_H^1$, ${\frak g}_H^2$ are responsible
for the existence of nonlinear automorphisms of $Q_H$ that preserve
the origin.

An example of how the algebra ${\frak g}_H$ is used in $CR$-geometry is
the equivalence problem for {\it strongly uniform} $CR$-manifolds. Let
$H_1,H_2$ be two $\RR^k$-valued Hermitian forms on $\CC^n\times\CC^n$. We say
that $H_1$ and $H_2$ are {\it equivalent}, if there exist linear
transformations $A$ of $\CC^n$ and $B$ of $\RR^k$ such that
$$
H_2(z,z)=BH_1(Az,Az).
$$
We call a $CR$-manifold $M$ strongly uniform,
if the forms ${\cal L}_M(p)$ are equivalent for all $p\in M$. If, for
example, $M$ is Levi non-degenerate and $CR\hbox{codim}M=1$ then $M$
is strongly uniform. The equivalence problem for strongly uniform Levi
non-degenerate integrable $CR$-manifolds is usually approached by
constructing a $CR$-invariant parallelism on certain bundles over the manifolds with
values in a suitable Lie algebra. In the number of cases (see \cite{C},
\cite{CM}, \cite{L}, \cite{EIS}) this Lie algebra was chosen to be
${\frak g}_H$, where $H$ is a Hermitian form equivalent to any
${\cal L}_M(p)$, $p\in M$. In the general approach of Tanaka \cite{Ta}, however, a
seemingly different algebra was used: Tanaka considered a certain maximal
prolongation $\tilde{\frak g}_H$ of ${\frak g}_H^{-2}\oplus{\frak
g}_H^{-1}\oplus{\frak g}_H^0$. It
is therefore a reasonable question whether the algebras ${\frak g}_H$ and
$\tilde{\frak g}_H$ are isomorphic. In this paper we give a positive
answer to this question in the main theorem below (see \cite {Ta}, \cite{L},
\cite{EIS} for partial results).

We will now give the precise definition of the algebra $\tilde {\frak
g}_H$ from \cite{Ta}. It is defined as an a priori
infinite-dimensional graded Lie algebra
$$
\tilde{\frak g}_H={\frak g}_H^{-2}\oplus{\frak g}_H^{-1}\oplus{\frak
g}_H^0\oplus
\left(\oplus_{l=1}^{\infty}\tilde{\frak g}_H^l\right)
$$
which is maximal among all Lie algebras of the above form that
satisfy the conditions:
\begin{eqnarray*}
&{}&\hbox{{\bf (i)} For $l\ge 0$ and $X\in\tilde{\frak g}_H^l$,
$[X,{\frak g}_H^{-1}]=0$ implies $X=0$};\\
&{}&\hbox{{\bf (ii)} ${\frak g}_H^{-2}\oplus{\frak g}_H^{-1}\oplus{\frak
g}_H^0$ is a subalgebra of $\tilde{\frak g}_H$.}
\end{eqnarray*}

It is shown in \cite{Ta} that $\tilde{\frak g}_H$ is unique and can be constructed
by the following inductive procedure. First we define vector
spaces $\tilde{\frak g}_H^l$ and brackets $[X_l,X_{-1}]\in\tilde{\frak
g}_H^{l-1}$, $[X_l,X_{-2}]\in\tilde{\frak
g}_H^{l-2}$, where $X_p\in\tilde{\frak g}_H^p$ (we set $\tilde{\frak
g}_H^l:={\frak g}_H^l$ for $l=-2,-1,0$). Suppose that these spaces and brackets have been defined
for for $0\le l\le L-1$ in such a way that the following holds
\begin{eqnarray*}
\left[\left[X_l,X_{-1}\right], Y_{-1}\right]-\left[\left[X_l,Y_{-1}\right],X_{-1}\right] &=& \left[X_l,\left[X_{-1},Y_{-1}\right]\right],\,\,\qquad\qquad\qquad(0.5.a)\\
\left[\left[X_l,X_{-2}\right],X_{-1}\right] &=& \left[\left[X_l,X_{-1}\right],X_{-2}\right],\qquad\qquad\qquad
(0.5.b)
\end{eqnarray*}
for all $X_l\in\tilde{\frak g}_H^l$, $X_{-1},Y_{-1}\in\tilde{\frak
g}_H^{-1}$. Then we define $\tilde{\frak g}_H^L$ to be the
vector space of all linear mappings $X_L:\tilde{\frak
g}_H^{-1}\rightarrow \tilde{\frak g}_H^{L-1}$ for which there exist
linear mappings $X'_L:\tilde{\frak g}_H^{-2}\rightarrow {\frak g}_H^{L-2}$
such that
\begin{eqnarray*}
\left[X_L(X_{-1}), Y_{-1}\right]-\left[X_L(Y_{-1}),
X_{-1}\right]&=&X_L'(\left[X_{-1},Y_{-1}\right]),\,\,\qquad\qquad\qquad
(0.6.a)\\
\left[X_L'(X_{-2}),X_{-1}\right]&=&\left[X_L(X_{-1}),X_{-2}\right],\qquad\qquad\qquad
(0.6.b)
\end{eqnarray*}
for all $X_{-1},Y_{-1}\in\tilde{\frak g}_H^{-1}$,
$X_{-2}\in\tilde{\frak g}_H^{-2}$. We set $[X_L,X_{-1}]:=X_L(X_{-1})$
for all $X_{-1}\in\tilde{\frak g}_H^{-1}$. Since $H$ is
non-degenerate, we have $\tilde{\frak g}_H^{-2}=[\tilde{\frak
g}_H^{-1},\tilde{\frak g}_H^{-1}]$, and therefore $X'_L$ is uniquely
determined by $X_L$. Then we set $[X_L,X_{-2}]:=X'_L(X_{-2})$ for all
$X_{-2}\in\tilde{\frak g}_H^{-2}$. We also set
$[X_{-1},X_l]:=-[X_l,X_{-1}]$ and $[X_{-2},X_l]:=-[X_l,X_{-2}]$.
Clearly, (0.6) then gives equations
(0.5) for $l=L$. 

Note that equations (0.5) imply
$$
[[X_l,X_{-2}],Y_{-2}]=[[X_l,Y_{-2}],X_{-2}],\eqno{(0.7)}
$$
for all $X_l\in\tilde{\frak g}_H^l$, $l\ge 0$, and
$X_{-2},Y_{-2}\in\tilde{\frak g}_H^{-2}$.

Let us now define brackets $[X_p,X_q]\in\tilde{\frak g}_H^{p+q}$,
$X_p\in\tilde{\frak g}_H^p$, $X_q\in\tilde{\frak g}_H^q$, $p,q\ge 0$,
inductively as follows. Suppose that these brackets have been defined
for $p,q\ge 0$, $p+q\le L-1$, in such a way that for any
$X_p\in\tilde{\frak g}_H^p$, $X_q\in\tilde{\frak g}_H^q$ the following holds
\begin{eqnarray*}
\left[\left[X_p,X_q\right],X_{-1}\right]&=&\left[\left[X_p,X_{-1}\right],X_q\right]+\left[X_p,\left[X_q,X_{-1}\right]\right],\qquad\qquad\qquad
(0.8.a)\\
\left[\left[X_p,X_q\right],X_{-2}\right]&=&\left[\left[X_p,X_{-2}\right],X_q\right]+\left[X_p,\left[X_q,X_{-2}\right]\right],\qquad\qquad\qquad
(0.8.b)
\end{eqnarray*}
for all $X_{-1}\in\tilde{\frak g}_H^{-1}$, $X_{-2}\in\tilde{\frak
g}_H^{-2}$. We take any $X_p\in\tilde{\frak g}_H^p$,
$X_q\in\tilde{\frak g}_H^q$ with $p,q\ge 0$ and $p+q=L$ and define
linear mappings $X_L$ and $X'_L$ from $\tilde{\frak g}_H^{-1}$ and $\tilde{\frak
g}_H^{-2}$ to ${\frak g}_H^{L-1}$ and ${\frak g}_H^{L-2}$ respectively
by
\begin{eqnarray*}
X_L(X_{-1})&:=&[[X_p,X_{-1}],X_q]+[X_p,[X_q,X_{-1}]],\\
X'_L(X_{-2})&:=&[[X_p,X_{-2}],X_q]+[X_p,[X_q,X_{-2}]].
\end{eqnarray*}
Then we see that $X_L,X'_L$ so defined satisfy (0.6) and therefore
$X_L\in\tilde{\frak g}_H^L$. We then define $[X_p,X_q]:=X_L$. Clearly,
this definition gives identities (0.8) for all $p,q\ge 0$,
$p+q=L$. Thus $[X_p,X_q]$ have been defined for all $p,q\ge 0$. Note that $[X_p,X_q]=-[X_q,X_p]$ for all $p,q\ge 0$. By
induction, we can also prove
$$
[[X_p,X_q],X_r]+[[X_q,X_r],X_p]+[[X_r,X_p],X_q]=0,\eqno{(0.9)}
$$
for all $X_p\in\tilde{\frak g}_H^p$, $X_q\in\tilde{\frak g}_H^q$,
$X_r\in\tilde{\frak g}_H^r$, $p,q,r\ge 0$. By (0.5), (0.7), (0.8),
(0.9) the brackets defined above give a Lie algebra structure on
$\tilde{\frak g}_H$. This completes the
construction of $\tilde{\frak g}_H$ in \cite{Ta}.

We now define a mapping $\Phi:{\frak g}_H\rightarrow\tilde{\frak
g}_H$ as follows:
\begin{eqnarray*}
&{}&\hbox{$\Phi$ is identical on ${\frak g}_H^{-2}\oplus{\frak
g}_H^{-1}\oplus
{\frak g}_H^0$},\\
&{}&\hbox{$[\Phi(X)](X_{-1}):=[X,X_{-1}]$ for $X\in{\frak g}_H^1$},\\
&{}&\hbox{$[[\Phi(X)](X_{-1})](Y_{-1}):=[[X,X_{-1}],Y_{-1}]$ for
$X\in{\frak
g}_H^2$}.
\end{eqnarray*}
It follows that $\Phi$ is a Lie algebra homomorphism and $\hbox{ker}\,\Phi=\{0\}$. Moreover, $\Phi({\frak g}_H^p)\subset
\tilde{\frak g}_H^p$ for $p=1,2$.

We are now ready to formulate the main result of the paper.

\begin{theorem} The mapping $\Phi$ is an isomorphism.
\end{theorem}

We will prove the theorem in the next section. Before proceeding, we
would like to thank G. Schmalz for useful discussions.

\section {Proof of Theorem}

It is clear from the preceding discussion that to prove the
theorem it is sufficient to show that $\hbox{dim}\,{\frak g}_H^p=
\hbox{dim}\,\tilde{\frak g}_H^p$ for $p=1,2$, and $\tilde{\frak
g}_H^3=\{0\}$.

\begin{lemma} $\hbox{dim}\,{\frak g}_H^1=
\hbox{dim}\,\tilde{\frak g}_H^1$.
\end{lemma}

{\bf Proof:} Let $X_1\in\tilde{\frak g}_H^1$. Since ${\frak
g}_H^{-1}$, ${\frak g}_H^0$
are given in the form as in (0.3.b), (0.3.c), $X_1$ can be written as
$$
X_1\left(p\frac{\partial}{\partial z}+2iH(z,p)\frac{\partial}{\partial
w}\right)=\phi(p)z\frac{\partial}{\partial
z}+\psi(p)w\frac{\partial}{\partial w},
$$
$p\in\CC^n$, where $\phi$, $\psi$ are real-linear mappings from
$\CC^n$ to the spaces $M(n,\CC)$ of complex $n\times n$- and
$M(k,\RR)$ of real $k\times k$-matrices 
respectively such that, for any $p,z\in\CC^n$,
$$
\psi(p)H(z,z)=2\hbox{Re}\,H(\phi(p)z,z) \eqno{(1.1)}
$$
(see (0.2.a)).
Let $X'_1$ be the linear mapping from ${\frak g}_H^{-2}$ to ${\frak
g}_H^{-1}$ corresponding to $X_1$ as in the definition of
$\tilde{\frak g}_H^1$. It then follows from (0.3.a), (0.3.b) that
$X'_1$ can be written in the form
$$
X'_1\left(q\frac{\partial}{\partial
w}\right)=\mu(q)\frac{\partial}{\partial
z}+2iH(z,\mu(q))\frac{\partial}{\partial w},
$$
$q\in\RR^k$, where $\mu$ is a linear mapping from $\RR^k$ to $\CC^n$.
Next, conditions (0.6) for $L=1$ are equivalent to
\begin{eqnarray*}
4\mu\left(\hbox{Im}\,H(p_1,p_2)\right)&=&\phi(p_1)p_2-\phi(p_2)p_1,\qquad\qquad
(1.2.a)\\
4\hbox{Im}\,H(\mu(q),p)&=&\psi(p)q,\,\,\,\,\,\,\,\qquad\qquad\qquad\qquad
(1.2.b)
\end{eqnarray*}
for all $p,p_1,p_2\in\CC^n$, $q\in\RR^k$.

We set
$$
A(p,p):=\frac{1}{2}\phi(p)p-i\mu(H(p,p)), \qquad a:=\mu,
$$
$p\in\CC^n$. We will show that the following holds (cf. (0.2.b))
$$
H(A(p,p),p)=2iH(p,aH(p,p)), \eqno{(1.3)}
$$
for all $p\in\CC^n$. We write $\phi$ in the most general form
$$
\phi(p)=Rp+Q\overline{p},
$$
where $R$, $Q$ are constant vectors of length $n$ with entries from
$M(n,\CC)$. Formulas (1.1), (1.2) then give
\begin{eqnarray*}
Rp_1p_2&=&Rp_2p_1,\qquad\qquad\qquad\qquad\qquad\qquad\qquad\qquad (1.4.a)\\
Q\overline{p_1}p_2&=&2i\mu(H(p_2,p_1)),\,\,\qquad\qquad\qquad\qquad\qquad\qquad (1.4.b)\\
H(Rp_1p_2,p_3)&=&2iH(p_1,\mu(H(p_2,p_3)))-H(p_2,Q\overline{p_1}p_3),\,\,\,\qquad(1.4.c)
\end{eqnarray*}
for all $p_1,p_2,p_3\in\CC^n$, where $\mu$ is complex-linearly extended from $\RR^k$ to
$\CC^k$. 
Identities (1.3) easily follow from (1.4.b), (1.4.c). Identity (1.4.b) in
addition gives
$$
A(p,p)=\frac{1}{2}Rpp,
$$
thus showing that $A$ is a quadratic form on $\CC^n\times\CC^n$.

It is clear from (1.1), (1.4.b), (1.4.c) that $a$ uniquely determines $X_1$ (also note
that (1.3) implies that $A$ and $a$ uniquely determine each other), and
the lemma is
proved.\hfill $\Box$

\begin{lemma} $\hbox{dim}\,{\frak g}_H^2=
\hbox{dim}\,\tilde{\frak g}_H^2$.
\end{lemma}

{\bf Proof:} Let $X_2\in\tilde{\frak g}_H^2$. It follows from 
(0.3.a)--(0.3.c) that
there exist real-bilinear mappings $\phi(\cdot,\cdot)$ and
$\psi(\cdot,\cdot)$ from $\CC^n\times\CC^n$ to the spaces $M(n,\CC)$ and $M(k,\RR)$ respectively, and a real-bilinear mapping
$\mu(\cdot,\cdot)$ from $\CC^n\times\RR^k$ to $\CC^n$ such that
\begin{eqnarray*}
X_2\left(p_1\frac{\partial}{\partial z}+2iH(z,p_1)\frac{\partial}{\partial
w}\right)
\left(p_2\frac{\partial}{\partial z}+2iH(z,p_2)\frac{\partial}{\partial
w}\right)
&=&\phi(p_1,p_2)z\frac{\partial}{\partial
z}+\psi(p_1,p_2)w\frac{\partial}{\partial w},\\
\Biggl[X_2\left(p\frac{\partial}{\partial z}+2iH(z,p)\frac{\partial}{\partial
w}\right)\Biggr]'
\left(q\frac{\partial}{\partial w}\right)
&=&\mu(p,q)\frac{\partial}{\partial
z}+2iH(z,\mu(p,q))\frac{\partial}{\partial w},\\
\end{eqnarray*}
$p,p_1,p_2\in\CC^n$, $q\in\RR^k$, where $\Biggl[X_2\left(p\frac{\partial}{\partial z}+2iH(z,p)\frac{\partial}{\partial
w}\right)\Biggr]'$ corresponds to $X_2\left(p\frac{\partial}{\partial z}+2iH(z,p)\frac{\partial}{\partial
w}\right)$ as an element of $\tilde{\frak g}_H^1$. Let $X'_2$ be the corresponding linear mapping from
${\frak g}_H^{-2}$ to ${\frak g}_H^0$. It follows from (0.3.a), (0.3.c) that it
can be written in the form
$$
X'_2\left(q\frac{\partial}{\partial
w}\right)=\eta(q)z\frac{\partial}{\partial
z}+\nu(q)w\frac{\partial}{\partial w},
$$
$q\in\RR^k$, where $\eta$ and $\nu$ are linear mappings from $\RR^k$
to the spaces $M(n,\CC)$ and $M(k,\RR)$
respectively. Equation (0.2.a) gives that the following conditions are
satisfied
\begin{eqnarray*}
\psi(p_1,p_2)H(z,z)&=&2\hbox{Re}\,H(\phi(p_1,p_2)z,z),\,\qquad\qquad\qquad
(1.5.a)\\
\nu(q)H(z,z)&=&2\hbox{Re}\,H(\eta(q)z,z),\qquad\qquad\qquad\qquad (1.5.b)
\end{eqnarray*}
for all $p_1,p_2,z\in\CC^n$, $q\in\RR^k$. Next, analogously to (1.2),
the following holds
\begin{eqnarray*}
4\mu\left(p_1,\hbox{Im}\,H(p_2,p_3)\right)&=&\phi(p_1,p_2)p_3-\phi(p_1,p_3)p_2,\qquad\qquad
(1.6.a)\\
4\hbox{Im}\,H(\mu(p_1,q),p_2)&=&\psi(p_1,p_2)q,\,\qquad\qquad\qquad\qquad\qquad
(1.6.b)
\end{eqnarray*}
for all $p,p_1,p_2\in\CC^n$, $q\in\RR^k$. 
Further, conditions (0.6) for
$L=2$ are equivalent to
\begin{eqnarray*}
4\eta\left(\hbox{Im}\,H(p_1,p_2)\right)&=&\phi(p_2,p_1)-\phi(p_1,p_2),\,\,\,\qquad\qquad
(1.7.a)\\
\mu(p,q)&=&-\eta(q)p,\qquad\qquad\qquad\qquad\qquad (1.7.b)
\end{eqnarray*}
for all $p,p_1,p_2\in\CC^n$, $q\in\RR^k$.

We set
$$
B(p,s):=\eta(s)p,\qquad r(s_1,s_2):=\frac{1}{2}\nu(s_2)s_1,
$$
$p\in\CC^n$,
$s,s_1,s_2\in\CC^k$, where $\eta$, $\nu$ are complex-linearly extended to
$\CC^k$. Then (1.5.b) implies
$$
\hbox{Re}\,H(B(p,q),p)=r(H(p,p),q), \eqno{(1.8)}
$$
for all $p\in\CC^n, q\in\RR^k$, which is analogous to (0.2.c). It follows from
(1.6) that $\phi$ is uniquely determined by $\mu$ (as in
(1.4) above). Therefore, by (1.5.a)
and (1.7.b), $X_2$ is uniquely
determined by $B$ (note that $B$ also uniquely determines $r$
by (1.8)). Thus, it is clear from (0.2.d) that to prove the lemma, we
need to show that
$$
\hbox{Im}\,H(B(p,H(p,p)),p)=0, \eqno{(1.9)}
$$
for all $p\in\CC^n$ and that $r(s_1,s_2)$ is symmetric. We write $\phi$ in the most general form
$$
\phi(p_1,p_2)=Mp_1p_2+N\overline{p_1}p_2+Tp_1\overline{p_2}+S\overline{p_1}\overline{p_2}.
$$
Then (1.5.a), (1.6), (1.7) give
\begin{eqnarray*}
M&=&0,\qquad\qquad\qquad\qquad\qquad\qquad\qquad\qquad\qquad (1.10.a)\\
S&=&0,\qquad\qquad\qquad\qquad\qquad\qquad\qquad\qquad\qquad (1.10.b)\\
N\overline{p_1}p_2p_3&=&N\overline{p_1}p_3p_2,\,\,\,\,\,\qquad\qquad\qquad\qquad\qquad\qquad\qquad
(1.10.c)\\
N\overline{p_1}p_2p_3-Tp_2\overline{p_1}p_3&=&2i\mu(p_3,H(p_2,p_1)),\,\,\,\,\,\,\qquad\qquad\qquad\qquad
\qquad(1.10.d)\\
Tp_1\overline{p_2}p_3&=&2i\mu(p_1,H(p_3,p_2)),\,\,\,\,\,\,\qquad\qquad\qquad\qquad\qquad (1.10.e)\\
H(N\overline{p_1}p_2p_3,p_4)&=&2iH(p_2,\mu(p_1,H(p_3,p_4)))-H(p_3,Tp_1\overline{p_2}p_4),\,\,
(1.10.f)
\end{eqnarray*}
for all $p_1,p_2,p_3,p_4\in\CC^n$, where $\mu$ is extended in the last
argument to a complex-linear mapping on $\CC^k$. Calculating
$\hbox{Im}\,H(\eta(H(p,p))p,p)$ from (1.7.a) we get
$$
\hbox{Im}\,H(\eta(H(p,p))p,p)=\frac{1}{2}\hbox{Re}\,H(N\overline{p}pp-Tp\overline{p}p,p).
$$
On the other hand, (1.5.a), (1.6.b), (1.7.b) give
$$
\hbox{Im}\,H(\eta(H(p,p))p,p)=-\frac{1}{2}\hbox{Re}\,H(N\overline{p}pp+Tp\overline{p}p+Mppp+S\overline{p}\overline{p}p,p).
$$
Comparing the last two expressions and using (1.10.a), (1.10.b), (1.10.d), (1.10.e) yields
(1.9).

To show that $r(s_1,s_2)$ is symmetric, by (1.5.b), we need to prove that
$$
\hbox{Re}\,H(\eta(H(p_1,p_1))p_2,p_2)=\hbox{Re}\,H(\eta(H(p_2,p_2))p_1,p_1),\eqno{(1.11)}
$$
for all $p_1,p_2\in\CC^n$. It follows from (1.7.a) that
$$
\hbox{Re}\,H(\eta(H(p_1,p_1))p_2,p_2)=-\frac{1}{2}\hbox{Im}\,H(N\overline{p_1}p_1p_2-Tp_1\overline{p_1}p_2,p_2),\eqno{(1.12)}
$$
for all $p_1,p_2\in\CC^n$. On the other hand, (1.5.a), (1.6.b),
(1.7.b), (1.10.a), (1.10.b) give
$$
\hbox{Re}\,H(\eta(H(p_1,p_1))p_2,p_2)=-\frac{1}{2}\hbox{Im}\,H(N\overline{p_2}p_2p_1-Tp_2\overline{p_2}p_1,p_1),
$$
for all $p_1,p_2\in\CC^n$, which together with (1.12) implies (1.11).

The lemma is proved.\hfill $\Box$

\begin{lemma} $\tilde{\frak g}_H^3=\{0\}$.
\end{lemma}

{\bf Proof:} Let $X_3\in\tilde{\frak g}_H^3$. Then there exist
real-trilinear
mappings $\phi(\cdot,\cdot,\cdot)$, $\psi(\cdot,\cdot,\cdot)$ from
$\CC^n\times\CC^n\times\CC^n$ to the spaces $M(n,\CC)$ and $M(k,\RR)$ respectively, real-bilinear mappings $\eta(\cdot,\cdot)$,
$\nu(\cdot,\cdot)$ from $\CC^n\times\RR^k$ to the above spaces of
matrices respectively, and a real-trilinear mapping
$\mu(\cdot,\cdot,\cdot)$ from $\CC^n\times\CC^n\times\RR^k$ to $\CC^n$
such that
\begin{eqnarray*}
&{}&X_3\left(p_1\frac{\partial}{\partial z}+2iH(z,p_1)\frac{\partial}{\partial
w}\right)
\left(p_2\frac{\partial}{\partial z}+2iH(z,p_2)\frac{\partial}{\partial
w}\right)\left(p_3\frac{\partial}{\partial z}+2iH(z,p_3)\frac{\partial}{\partial
w}\right)=\\
&{}&\phi(p_1,p_2,p_3)z\frac{\partial}{\partial
z}+\psi(p_1,p_2,p_3)w\frac{\partial}{\partial w},\\
&{}&\Biggl[X_3\left(p_1\frac{\partial}{\partial z}+
2iH(z,p_1)\frac{\partial}{\partial
w}\right)\left(p_2\frac{\partial}{\partial z}+2iH(z,p_2)\frac{\partial}{\partial
w}\right)
\Biggr]'
\left(q\frac{\partial}{\partial w}\right)=\\
&{}&\mu(p_1,p_2,q)\frac{\partial}{\partial
z}+2iH(z,\mu(p_1,p_2,q))\frac{\partial}{\partial w},\\
&{}&\Biggl[X_3\left(p\frac{\partial}{\partial z}+2iH(z,p)\frac{\partial}{\partial
w}\right)\Biggr]'
\left(q\frac{\partial}{\partial w}\right)=\\
&{}&\eta(p,q)z\frac{\partial}{\partial
z}+\nu(p,q)w\frac{\partial}{\partial w},
\end{eqnarray*}
$p,p_1,p_2,p_3\in\CC^n$, $q\in\RR^k$, where $[X_3\left(p_1\frac{\partial}{\partial z}+
2iH(z,p_1)\frac{\partial}{\partial
w}\right)\left(p_2\frac{\partial}{\partial z}+2iH(z,p_2)\frac{\partial}{\partial
w}\right)
\Biggr]'$ corresponds to $X_3\left(p_1\frac{\partial}{\partial z}+
2iH(z,p_1)\frac{\partial}{\partial
w}\right)\left(p_2\frac{\partial}{\partial z}+2iH(z,p_2)\frac{\partial}{\partial
w}\right)$ as an element of $\tilde{\frak g}_H^1$ and $[X_3\left(p\frac{\partial}{\partial z}+2iH(z,p)\frac{\partial}{\partial
w}\right)\Biggr]'$ corresponds to $X_3\left(p\frac{\partial}{\partial z}+2iH(z,p)\frac{\partial}{\partial
w}\right)$ as an element of $\tilde{\frak g}_H^2$. Let $X'_3$ be
the corresponding linear mapping from ${\frak g}_H^{-2}$ to
$\tilde{\frak g}_H^1$. Then there exist real-bilinear mappings
$\lambda(\cdot,\cdot)$ and $\rho(\cdot,\cdot)$ from $\CC^n\times\RR^k$
to the spaces $M(n,\CC)$ and $M(k,\RR)$ respectively such that
$$
X'_3\left(q\frac{\partial}{\partial
w}\right)\left(p\frac{\partial}{\partial
z}+2iH(z,p)\frac{\partial}{\partial w}\right)=
\lambda(p,q)z\frac{\partial}{\partial
z}+\rho(p,q)w\frac{\partial}{\partial w},
$$
$p\in\CC^n$, $q\in\RR^k$. Equation (0.2.a) gives
\begin{eqnarray*}
\psi(p_1,p_2,p_3)H(z,z)&=&2\hbox{Re}\,H(\phi(p_1,p_2,p_3)z,z),\qquad\qquad\qquad
(1.13.a)\\
\nu(p,q)H(z,z)&=&2\hbox{Re}\,H(\eta(p,q)z,z),\,\qquad\qquad\qquad\qquad (1.13.b)\\
\rho(p,q)H(z,z)&=&2\hbox{Re}\,H(\lambda(p,q)z,z),\,\qquad\qquad\qquad\qquad
(1.13.c)
\end{eqnarray*}
for all $p,p_1,p_2,p_3,z\in\CC^n$, $q\in\RR^k$. Next, analogously to
(1.6), we have
\begin{eqnarray*}
4\mu\left(p_1,p_2,\hbox{Im}\,H(p_3,p_4)\right)&=&\phi(p_1,p_2,p_3)p_4-\phi(p_1,p_2,p_4)p_3,\,\,\,\qquad
(1.14.a)\\
4\hbox{Im}\,H(\mu(p_1,p_2,q),p_3)&=&\psi(p_1,p_2,p_3)q,\qquad\qquad\qquad\qquad\qquad
(1.14.b)
\end{eqnarray*}
for all $p_1,p_2,p_3\in\CC^n$, $q\in\RR^k$. Further, there are the
following analogues of identities (1.7)
\begin{eqnarray*}
4\eta\left(p_1,\hbox{Im}\,H(p_2,p_3)\right)&=&\phi(p_1,p_3,p_2)-\phi(p_1,p_2,p_3),\qquad
(1.15.a)\\
\mu(p_1,p_2,q)&=&-\eta(p_1,q)p_2,\,\,\qquad\qquad\qquad\qquad (1.15.b)
\end{eqnarray*}
for all $p_1,p_2\in\CC^n$, $q\in\RR^k$. Finally, conditions (0.6) for
$L=3$ are equivalent to
\begin{eqnarray*}
4\lambda(p_1,\hbox{Im}\,H(p_2,p_3))&=&\phi(p_3,p_2,p_1)-\phi(p_2,p_3,p_1),\qquad
(1.16.a)\\
\lambda(p,q)&=&\eta(p,q),\,\,\,\qquad\qquad\qquad\qquad\qquad (1.16.b)
\end{eqnarray*}
for all $p,p_1,p_2,p_3\in\CC^n$, $q\in\RR^k$. 

We will now show that identities (1.13)--(1.16) imply that $X_3=0$. It
follows by the argument in Lemma 1.2 from (1.13.a) and identities
(1.14), (1.15) that $\phi$ has the form (see (1.10.a), (1.10.b)) 
$$
\phi(p_1,p_2,p_3)=Gp_1\overline{p_2}p_3+K\overline{p_1}\overline{p_2}p_3+Fp_1p_2\overline{p_3}+L\overline{p_1}p_2\overline{p_3},
$$
and (see (1.10.c)--(1.10.f))
\begin{eqnarray*}
Gp_1p_2p_3p_4&=&Gp_1p_2p_4p_3,\,\,\qquad\qquad\qquad\qquad\qquad\qquad
(1.17.a)\\
(Gp_1+K\overline{p_1})\overline{p_2}p_3p_4&-&(Fp_1+L\overline{p_1})p_3\overline{p_2}p_4=\\
&{}&2i\mu(p_1,p_4,H(p_3,p_2)),\qquad\qquad\qquad\qquad (1.17.b)\\
(Fp_1+L\overline{p_1})p_2\overline{p_3}p_4&=&2i\mu(p_1,p_2,H(p_4,p_3)),\qquad\qquad\qquad\qquad
(1.17.c)\\
H((Gp_1+K\overline{p_1})\overline{p_2}p_3p_4,p_5)&=&2iH(p_3,\mu(p_1,p_2,H(p_5,p_4)))-\\
&{}&H(p_4,(Fp_1+L\overline{p_1})p_2\overline{p_3}p_5),
\,\,\qquad\qquad\qquad
(1.17.d)\\
Kp_1p_2p_3p_4&=&Kp_1p_2p_4p_3,\,\,\qquad\qquad\qquad\qquad\qquad\qquad
(1.17.e)
\end{eqnarray*}
for all $p_1,p_2,p_3,p_4\in\CC^n$, where $\mu$ is extended
complex-linearly to $\CC^k$ in the last argument. Further,
(1.15.b), (1.16) imply
$$
L\overline{p_1}p_2\overline{p_3}p_4-Gp_2\overline{p_1}p_3p_4=2i\mu(p_3,p_4,H(p_2,p_1)),\eqno{(1.18)}
$$
for all $p_1,p_2,p_3,p_4,p_5\in\CC^n$. From (1.17.b)--(1.17.d), (1.18) we obtain
\begin{eqnarray*}
H(K\overline{p_1}\overline{p_2}p_3p_4,p_5)&=&H(p_3,Gp_5\overline{p_4}p_1p_2)-
H(p_4,Fp_1p_2\overline{p_3}p_5),\,\,\,\,\,\,
(1.19.a)\\
H(Gp_1\overline{p_2}p_3p_4,p_5)&=&-H(p_3,L\overline{p_4}p_5\overline{p_1}p_2)-
H(p_4,L\overline{p_1}p_2\overline{p_3}p_5),\,\,\,
(1.19.b)\\
Fp_1p_2p_3p_4&=&-Gp_4p_3p_1p_2,\,\,\,\,\,\,\,\,\,\,\,\,\qquad\qquad\qquad\qquad\qquad
(1.19.c)\\
Lp_1p_2p_3p_4&=&Lp_3p_4p_1p_2,\,\,\,\,\,\,\qquad\qquad\qquad\qquad\qquad\qquad
(1.19.d)\\
Kp_1p_2p_3p_4&=&Lp_1p_3p_2p_4+Lp_1p_4p_2p_3,\,\,\,\,\,\,\,\,\,\qquad\qquad\qquad
(1.19.e)\\
Gp_1p_2p_3p_4+Gp_3p_2p_1p_4&=&Fp_3p_1p_2p_4,\,\,\,\,\,\,\qquad\qquad\qquad\qquad\qquad\qquad
(1.19.f)
\end{eqnarray*}
for all $p_1,p_2,p_3,p_4,p_5\in\CC^n$.

We set
\begin{eqnarray*}
D(p_1,p_2,H(p_3,p_4))&:=&iGp_3\overline{p_4}p_2p_1,\,\,\,\,\,\,\qquad\qquad\qquad\qquad
(1.20.a)\\
t(H(p_1,p_2),H(p_3,p_4))&:=&-\frac{1}{2}L\overline{p_2}p_1\overline{p_4}p_3,\qquad\qquad\qquad\qquad
(1.20.b)
\end{eqnarray*}
for all $p_1,p_2,p_3,p_4\in\CC^n$. It follows from (1.17.a), (1.18),
(1.19.d) and the non-degeneracy of $H$ that (1.20.a) defines a
complex-trilinear form $D$ on $\CC^n\times\CC^n\times\CC^k$ symmetric
with respect to the first two variables and (1.20.b) defines a complex-bilinear
symmetric form $t$ on $\CC^k\times\CC^k$, both valued in $\CC^n$. It follows from (1.13.a), (1.19.c),
(1.19.e) that $X_3$ is uniquely determined by $t$ (note that it
follows from (1.19.b) that $D$ is uniquely determined by
$t$). The forms $D$ and $t$ satisfy the following relations
\begin{eqnarray*}
H(D(p,p,q),p)&=&4iH(p,t(q,H(p,p)),\,\qquad\qquad\qquad (1.21.a)\\
H(D(p,p,H(p,p)),p)&=&0,\qquad\qquad\qquad\qquad\qquad\qquad\qquad(1.21.b)
\end{eqnarray*}
for all $p\in\CC^n,q\in\RR^k$. Indeed, (1.21.a) follows from (1.19.b),
(1.19.d); to prove (1.21.b) we note that it follows from (1.19.c),
(1.19.f) that $Gp\overline{p}pp\equiv 0$.

We will now show that equations (1.21) can have only zero
solutions. For this we note that a polynomial vector field $X$ on $\CC^{n+k}$ of
the form
$$
X=\Bigl(D'(z,z,w)+t'(w,w)\Bigr)\frac{\partial}{\partial
z}+2iH\Bigl(z,t'(\overline{w},\overline{w})\Bigr)\frac{\partial}{\partial w}
$$
defines an infinitesimal holomorphic automorphism of $Q_H$ (i.e. $X\in{\frak g}_H$) if and only if the
following conditions are satisfied
\begin{eqnarray*}
H(D'(z,z,u),z)&=&4iH(z,t'(u,H(z,z)),\qquad\qquad\qquad (1.22.a)\\
H(D'(z,z,H(z,z)),z)&=&0,\,\qquad\qquad\qquad\qquad\qquad\qquad\qquad (1.22.b)
\end{eqnarray*}
for all $z\in\CC^n,u\in\RR^k$, where $D'$ is a complex-trilinear form
on $\CC^n\times\CC^n\times\CC^k$ symmetric with respect to the first
two variables, and $t'$ is a complex-bilinear symmetric form on
$\CC^k\times\CC^k$, both valued in $\CC^n$. Since the vector field $X$
has weight 3, it
must be zero by (0.1). This means that equations (1.22) can have only zero
solutions and therefore equations (1.21) have only zero solutions also.

Thus, $X_3=0$, and the lemma is proved.\hfill $\Box$

\bigskip

{\obeylines
Department of Pure Mathematics
The University of Adelaide
Box 498, G.P.O. Adelaide
South Australia 5001
AUSTRALIA
E-mail address: vezhov@spam.maths.adelaide.edu.au
\smallskip

Centre for Mathematics and Its Applications 
The Australian National University 
Canberra, ACT 0200
AUSTRALIA 
E-mail address: Alexander.Isaev@anu.edu.au
}

\end{document}